\documentclass{article}
\usepackage[utf8]{inputenc}
\usepackage{amssymb,graphicx,amsmath}
\usepackage{caption}
\usepackage[%
pdfauthor={A. Ibeas},%
pdftitle={Periodic lozenge tilings of the plane},%
hyperindex,plainpages=false,bookmarks=true,a4paper,%
pdfhighlight=/I,colorlinks=false,linkbordercolor={0.5 0.5 1},%
linkcolor=blue,citebordercolor={0.9 0.8 0.8},%
pdfborder={1 1 1},
]{hyperref}

\renewcommand\L{\Lambda}
\newcommand\Z{\mathbb Z}
\newcommand\N{\mathbb N}
\newcommand\R{\mathbb R}
\newcommand\D{\delta}
\renewcommand\Xi{\Delta}
\renewcommand\t{\tau}
\newcommand\ve[1]{{\mathbf #1}}

\newcommand\up{\(\bigtriangleup\) }
\newcommand\down{\(\bigtriangledown\) }

\newtheorem{defi}{Definition}

\newtheorem{teo}[defi]{Theorem}
\newtheorem{prop}[defi]{Proposition}

\newenvironment{demo}{\noindent{\bf Proof. }}{\hfill\(\blacksquare\)}

\begin{document}

\centerline{\large\bf Periodic lozenge tilings of the plane}\vspace{4mm}

\centerline{Álvar Ibeas Martín}
\centerline{Universidad de Cantabria}\vspace{1cm}

\begin{center}
\parbox{11cm}{{\bf Abstract.} This article addresses the problem of
  enumerating the tilings of a plane by lozenges, under the
  restriction that these tilings be doubly periodic. Kasteleyn's
  Pfaffian method is applied to compute the generating function of
  those permutations. The monomials of this function represent the
  different {\it types} of tilings, grouping them according to the
  number of lozenges in each orientation. We present an alternative
  approach to compute these types. Finally, two additional classes of
  tilings are proposed as open enumeration problems.}
\end{center}

\section{Introduction}

We consider the tiling of the plane by equilateral triangles, assuming
that the vertices of those are the points of the lattice \(\L_0\) spanned by
\(\ve u=(1,0)\), \(\ve v=(1/2,\sqrt 3/2)\).
\begin{figure}[h!]
\begin{center}
\includegraphics{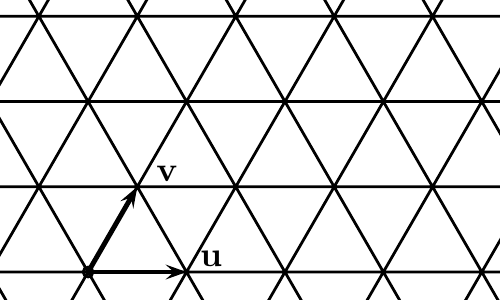}
\caption{Equilateral triangle tiling}\label{reja}
\end{center}
\end{figure}

\begin{figure}[h!]
\begin{center}
\includegraphics[width=3cm]{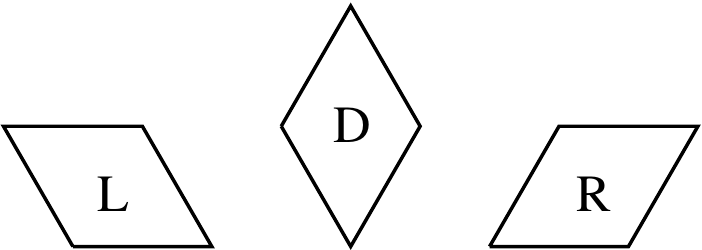}
\caption{}\label{notacion_2}
\end{center}
\end{figure}

By merging two adjacent triangles of this tiling we obtain a rhombus
or {\it lozenge}. It is obvious then that the plane can also be tiled
with such lozenges. These tiles may take three {\it orientations}, which
we denote by the notation set in Figure~\ref{notacion_2}.
\captionsetup{labelsep=none}
\captionsetup{labelsep=colon}

There is a rich literature devoted to studying the different lozenge
arrangements that tile certain bounded regions of the plane (see
\cite{Kup_94,Coh_Lar_Pro_98,Des_03,Kra_90} and the references
therein). The simplest scenario appears in the case of a convex
hexagon drawn over the grid of Figure~\ref{reja} with pairs of
opposite sides of the same length (a semiregular hexagon). Its lozenge
tilings are associated by an appealing bijection with the plane
partitions whose parts are bounded by the lengths of the hexagon
sides. We recall that the solid Young diagram of plane partition is an
arrangement of unit cubes located in the positive octant of \(\R^3\)
and satisfying that if there is a cube at position \((i,j,k)\),
\(a\leq i\), \(b\leq j\), and \(c\leq k\), then position \((a,b,c)\)
is neither empty.
\begin{figure}
\begin{center}
\includegraphics[height=3cm]{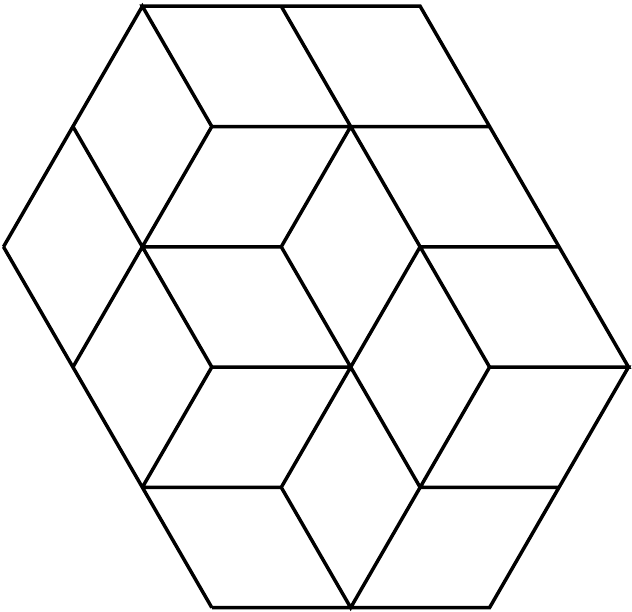}
\caption{A lozenge tiling in a semiregular hexagon}\label{fig_Y}
\end{center}
\end{figure}

A classical formula (proved in
\cite{Mac_60}) states that the number of plane partitions fitting in a
box of sides \(a,b\), and \(c\) is:
\[\prod_{i=0}^{a-1}\prod_{j=0}^{b-1}\prod_{k=0}^{c-1}\frac{i+j+k+2}{i+j+k+1}.\]

Throughout this article, we employ bold letters to denote vectors
\(\ve a=(a_1,a_2)\), considered as columns when using matrix notation:
\(\ve a=(a_1,a_2)^t\). We are interested in lozenge tilings of the
whole plane which are doubly periodic. More explicitly, let
\[B=\left[\begin{array}{c|c}
a_1&b_2\\ a_2&b_2\end{array}\right]\in\Z^{2\times 2}\]
be a 2-rank matrix and consider the sublattice
\begin{equation}\label{Lambda}
\L=[\ve u|\ve v]B\Z^2\subseteq\L_0.
\end{equation}
The tilings we are interested of are invariable under the translations
by \(a_1\ve u+a_2\ve v\) and \(b_1\ve u+b_2\ve v\); or
equivalently, by any element of \(\L\). Dealing with infinite tilings,
unlike in the finite regions case, we do not have a boundary to
``support'' our arguments.
\begin{figure}[h!]
\begin{center}
\includegraphics[width=5cm]{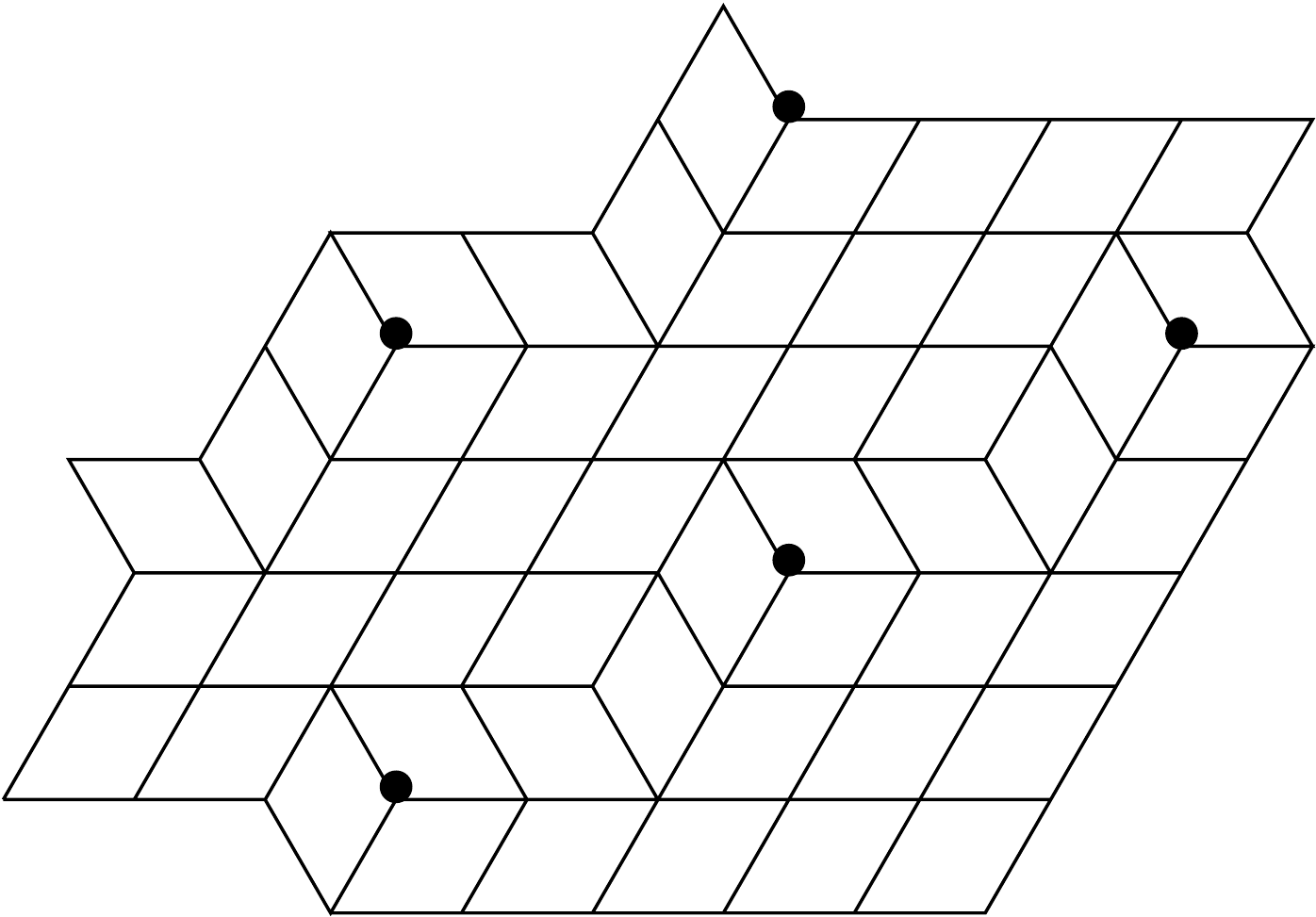}
\caption{A doubly periodic tiling}\label{mixta}
\end{center}
\end{figure}

We will define the {\it type} of a tiling as the number of lozenges in
each orientation it consists of. Note that, as can be checked using
the argument exposed in \cite{Dav_Tom_89}; in the mentioned case of a
semiregular hexagon with sides \(a\), \(b\), and \(c\), all its
lozenge tilings have a constant number of lozenges in each orientation
(\(ab\), \(ac\), and \(bc\)).

As we have said, a lozenge tiling is built merging pairs of triangles
under certain restrictions. Therefore, lozenge tilings correspond to
perfect matchings in an associated graph (see
Figure~\ref{eje_1}). Perfect matchings enumeration is a rich field of
research (see \cite{Pro_90} for an illustrating survey), stimulated by
several problems in the domain of Physics and Chemistry.

Kasteleyn developed a method for enumerating the perfect matchings of
a planar graph by means of a Pfaffian
computation~\cite{Kas_61,Kas_63,Kas_67}. In the first of these
references, the problem is also solved for a non-planar graph, which
can be embedded in a torus. This is accomplished by computing a linear
combination of four Pfaffians. In Section~\ref{sec_per_det}, we
explain how this methods applies to counting doubly periodic lozenge
tilings. The output is a polynomial {\it generating function}
\(Z(L,D,R)\in\Z[L,D,R]\), where the coefficient in \(L^iD^jR^k\)
equals the number of tilings of type \((i,j,k)\).

In Section~\ref{sec_tri}, we introduce a different approach to compute
the different types of doubly periodic tilings modulo a given period
\(\L\). We prove that the pair of height increments in the (infinite
analogue to the) solid Young diagram through a pair of vectors
spanning \(\L\) characterises the type of the tiling. As a
consequence, the possible types correspond to lattice points in a
certain triangle associated to a lattice basis \(B\), in the sense
of~(\ref{Lambda}). The vertices of this triangle correspond to the
three uniform tilings (those with all of its lozenges arranged in a
constant orientation).

Finally, in Section~\ref{sec_group}, we propose two ways of grouping
together significantly similar tilings, whose enumeration remains, up
to our knowledge, an open problem. Before all, let us fix the notation
and definitions we will need:

\section{Definitions}\label{sec_def}

There are two types of triangles in Figure~\ref{reja}: upwards and
downwards-pointing. Any lozenge contains a triangle of each type. In
this article, we refer to generic elements of these classes by the
symbols \up and \down\hspace{-1mm}. We identify both sets of triangles
with \(\L_0\) by means of the following convention: a point in
\(\L_0\) represents the right-most upwards-pointing and
downwards-pointing triangles which have that point as vertex.

Let \(\L\) be defined by Equation~(\ref{Lambda}). We refer to the
index of \(\L\) in \(\L_0\) simply as the {\it index} of \(\L\):
\[[\L_0:\L]=|\det B|=2(\mathrm{vol}\,\L)/\sqrt 3.\]
In order to define a \(\L\)-periodic tiling, we need to decide with
which of the three adjacent \down is merged every \up in a fundamental
set whose size is the index of \(\L\). The following map returns the
difference \down\hspace{-1mm}-\up within a lozenge on its
orientation as input.

\begin{center}\begin{tabular}{c@{\hspace{2cm}}c}\(\begin{array}{cccc}
\xi:&\{\mathrm{L,D,R}\}&\rightarrow&\L_0\\
&\mathrm L&\mapsto&\ve v-\ve u\\
&\mathrm D&\mapsto&\ve 0\\
&\mathrm R&\mapsto&\ve v.
\end{array}\)&\mbox{\parbox{5cm}{\vspace{-2mm} \includegraphics[scale=.5]{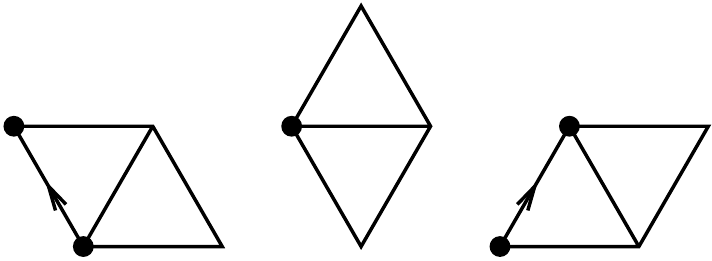}}}\end{tabular}
\end{center}

\renewcommand\theenumi{\roman{enumi}}
\renewcommand\labelenumi{\theenumi)}

\begin{figure}[h!]
\begin{center}
\includegraphics[width=2cm]{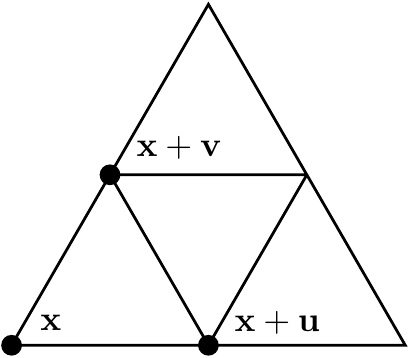}
\caption{Axiom {\it ii}') Exactly one of the three possible lozenges occurs.}\label{ii_1}
\end{center}
\end{figure}

\begin{defi}\label{def_tiling}
Let \(\L\subseteq\L_0\) be a 2-rank lattice. We define a {\bf
  \(\L\)-periodic tiling} as a map
\[\t:\L_0\rightarrow\{\mathrm{L,D,R}\}\]
satisfying the following two axioms:
\begin{enumerate}
\item {\bf Compatibility:}
\[\forall\ve x\in\L_0\,\forall\ve y\in\L,\ \t(\ve x)=\t(\ve x+\ve y).\]
\item {\bf Tiling:} The following map is bijective:
\[\begin{array}{cccc}
\tilde\t:&\L_0&\rightarrow&\L_0\\
&\ve x&\mapsto&\ve x+\xi\left(\t(\ve x)\right).\end{array}\]
\end{enumerate}
We use the notation \(T_\L\) for the set of \(\L\)-periodic tilings.
\end{defi}
Note that the second axiom is equivalent to (see Figure~\ref{ii_1}):
\renewcommand\labelenumi{\it ii')}
\begin{enumerate}
\item \(\forall\ve x\in\L_0\), exactly one of the following conditions is satisfied:
\[\tau(\ve x)=R,\ \tau(\ve x+\ve u)=L,\ \tau(\ve x+\ve v)=D.\]
\end{enumerate}
As we have already mentioned, it is enough to define \(\t\) on a
representative of each class modulo \(\L\). We consider the induced
mapping:
\[\begin{array}{ccccc}
\hat\t:&\L_0/\L&\rightarrow&\{\mathrm{L,D,R}\}\end{array}\]
and define the {\it type} of a \(\L\)-periodic tiling as:
\[t(\t,\L)=\left(\#\hat\t^{-1}(L),\#\hat\t^{-1}(D),\#\hat\t^{-1}(R)\right)\in\N^3.\]

\begin{figure}[h!]
\begin{center}
\includegraphics[width=6cm]{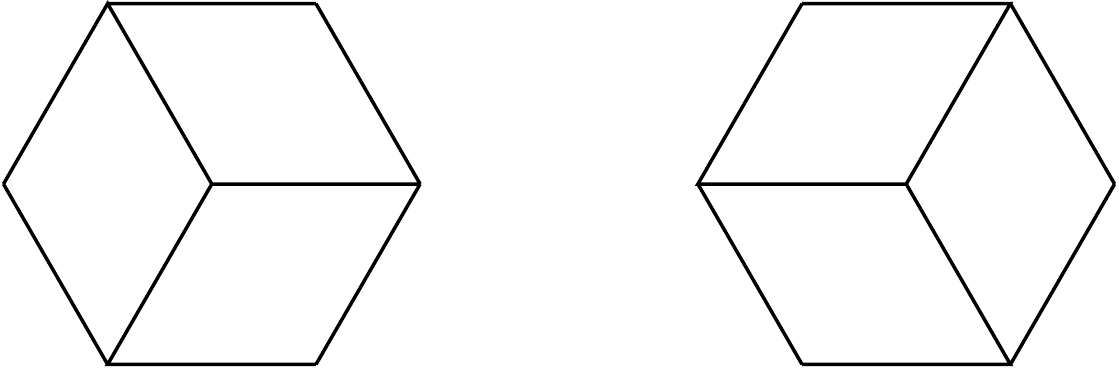}
\caption{An outside and an inside corner.}\label{esquinas}
\end{center}
\end{figure}

Infinite lozenge tilings can be represented by an infinite analogue to
a solid Young diagram, namely, the complement of a subset of \(\Z^3\)
closed under addition of elements in the semigroup \(\N^3\). In other
words, a staircase diagram in three dimensions; with two possible
identifications. As Figure~\ref{esquinas} shows, there are two lozenge
tilings of the unit hexagon. One of them is prominent (a solid cube),
and the other one is a ``hole'', limited by three walls. In the rest
of this article, we consider that the left design of
Figure~\ref{esquinas} is the solid cube.

With this identification, the vertices of a lozenge tiling can be
labelled by three-dimensional coordinates and given a {\it height
  function}, defined by the sum of these coordinates. The application
of height labels in this context dates back to
\cite{Blo_Hil_82,Thu_90} and has been useful in the study of tilings
of bounded regions.

An edge in the tiling is associated with a coordinates and a height
increment:
\[\begin{array}{ccc}\begin{array}{c|c|c|}
&F&f\\ \hline
\ve u&(-1,0,0)&-1\\
-\ve u&(+1,0,0)&+1\\
\ve v&(0,0,+1)&+1\\
-\ve v&(0,0,-1)&-1\\
\ve u-\ve v&(0,+1,0)&+1\\
\ve v-\ve u&(0,-1,0)&-1\\ \hline
\end{array}&\phantom{AAA}&\parbox{4cm}{
\includegraphics[width=3cm]{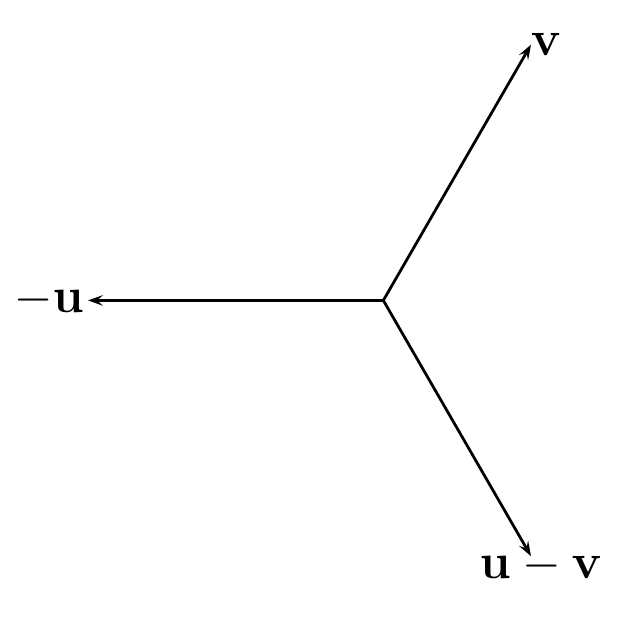}
}
\end{array}\]

Note that this effect of an edge on the height of a vertex is only
valid is the edge does appear in the lozenge tiling. For instance, if
two vertices \(\ve x\) and \(\ve x+\ve u\) are not connected by
\(\ve u\), then \(\t(\ve x)=D\) and they are connected by the
concatenation \((\ve v,\ve u-\ve v)\). Therefore, the coordinates
increment is not \((-1,0,0)\), but \((0,+1,+1)\).  In general, if the
height increment of an existing edge is \(h\), the effect of the same
non-occurring edge is \(-2h\).

\begin{defi}\label{def_path}
Let \(\L\subseteq\L_0\) be a 2-rank lattice and \(\t\) an
\(\L\)-periodic tiling. We define a {\bf path} in \(\t\) as a succession
\((\ve x_i)_{i=0}^N\subseteq\L_0\) satisfying
\[\ve
x_{i+1}-\ve x_i\in\{\pm\ve u,\pm\ve v,\pm(\ve u-\ve v)\}\] and:
\begin{itemize}
\item \(\ve x_{i+1}-\ve x_i=\ve u\Rightarrow\tau(\ve x_i)\neq D\)
\item \(\ve x_{i+1}-\ve x_i=\ve v\Rightarrow\tau(\ve x_i)\neq L\)
\item \(\ve x_{i+1}-\ve x_i=\ve v-\ve u\Rightarrow\tau(\ve x_i-\ve u)\neq R\)
\item \(\ve x_{i+1}-\ve x_i=-\ve u\Rightarrow\tau(\ve x_i-\ve u)\neq D\)
\item \(\ve x_{i+1}-\ve x_i=-\ve v\Rightarrow\tau(\ve x_i-\ve v)\neq L\)
\item \(\ve x_{i+1}-\ve x_i=\ve u-\ve v\Rightarrow\tau(\ve x_i-\ve v)\neq R\)
\end{itemize}
\end{defi}
Each path \(p=(\ve x_i)_{i=0}^{N}\) in \(\t\) can be associated with
a height increment, in the following way:
\[h_\t(p)=\sum_{i=0}^{N-1}f(\ve x_{i+1}-\ve x_i).\]

For a given tiling \(\t\), the height increment of a path only depends
on its extreme points. Therefore, setting \(h_\t(0)=0\), we can
associate a height to every point in \(\L_0\), defining a mapping
\(h_\t\) over \(\L_0\). As we usually treat elements of \(\L_0\) by
its coordinates on basis \((\ve u|\ve v)\), we will employ the
following change of coordinates:
\[\begin{array}{cccc}
e_\t:&\Z^2&\longrightarrow&\Z\\
&\ve a&\mapsto&h_\t([\ve u|\ve v]\ve a).
\end{array}\]
\begin{figure}[h!]
\begin{center}
\includegraphics[scale=.8]{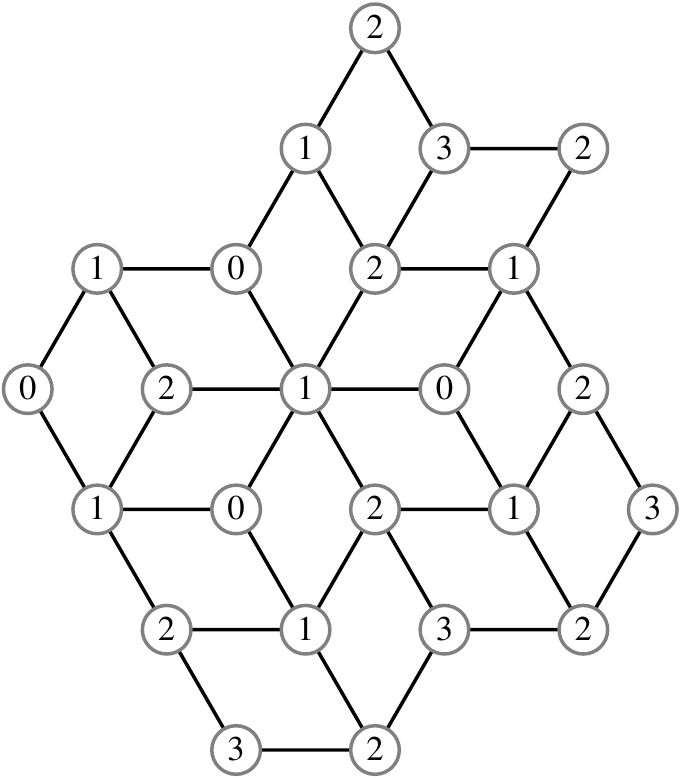}
\caption{The \(h\) function.}\label{alturas}
\end{center}
\end{figure}

Note that as a lozenge tiling is the projection of an pile of cubes
along lines parallel to vector \((1,1,1)\), the height \(h_\t(\ve x)\)
of a point (together with its projected position in the plane)
determines its three-dimensional coordinates, which by analogy we may
denote \(H_\t(\ve x)\).

\begin{prop}\label{esqueleto}Let \(\L\subseteq\L_0\) be a 2-rank lattice, \(\t\) a
  \(\L\)-periodic lattice and \(\ve x=\lambda_1\ve u+\lambda_2\ve
  v\in\L_0\). Then,
\[H_\t(\ve x)=(-\lambda_1,0,\lambda_2)+\frac{h_\t(\ve x)-\lambda_2+\lambda_1}{3}(1,1,1).\]
\end{prop}
\begin{demo} The maybe improper ``path'' \((0,\ve u,2\ve
  u,\ldots,\lambda_1\ve u,\lambda_1\ve u+\ve v,\ldots\ve x)\)
  would give a height label \(h=\lambda_2-\lambda_1\) and coordinates
  label \(H=(-\lambda_1,0,\lambda_2)\). Any wrong step can be replaced
  by a concatenation of two, changing the height by \(\pm 3\) and the
  coordinates by \(\pm(1,1,1)\).
\end{demo}

Note that if \(\t\) and \(\t'\) are tilings, we have \(H_\t(\ve
x)-H_{\t'}(\ve x)\in\Z\langle(1,1,1)\rangle\) and \(h_\t(\ve
x)-h_{\t'}(\ve x)\in(3)\), for every \(\ve x\in\L_0\).

Let \(B=[\ve a|\ve b]\) be a basis of \(\L\) with respect to
\(\L_0\) (see Equation~(\ref{Lambda})). We say that the pair
\(\D(\t,B)=(e_\t(\ve a),e_\t(\ve b))\) composed of the height
increment of the basic vectors is the {\it fingerprint} of \(\t\) in
basis \(B\). We note that for \(\ve x\in\L_0\) and \(\ve y\in\L\),
\[h_\tau(\ve x+\ve y)=h_\tau(\ve x)+h_\tau(\ve y),\]
and therefore,
\[h_\tau(\ve x+[\ve u|\ve v]B(\lambda_1,\lambda_2)^t)=h_\tau(\ve
x)+\lambda_1e_\tau(\ve a)+\lambda_2e_\tau(\ve
b),\ \forall\lambda_1,\lambda_2\in\Z.\]
In other words, the fingerprint of a tiling determines the height of
every point in the lattice \(\L\).

\section{The Permanent-Determinant method}\label{sec_per_det}

\begin{figure}
\begin{center}
\begin{tabular}{cc}
\includegraphics[height=2.7cm]{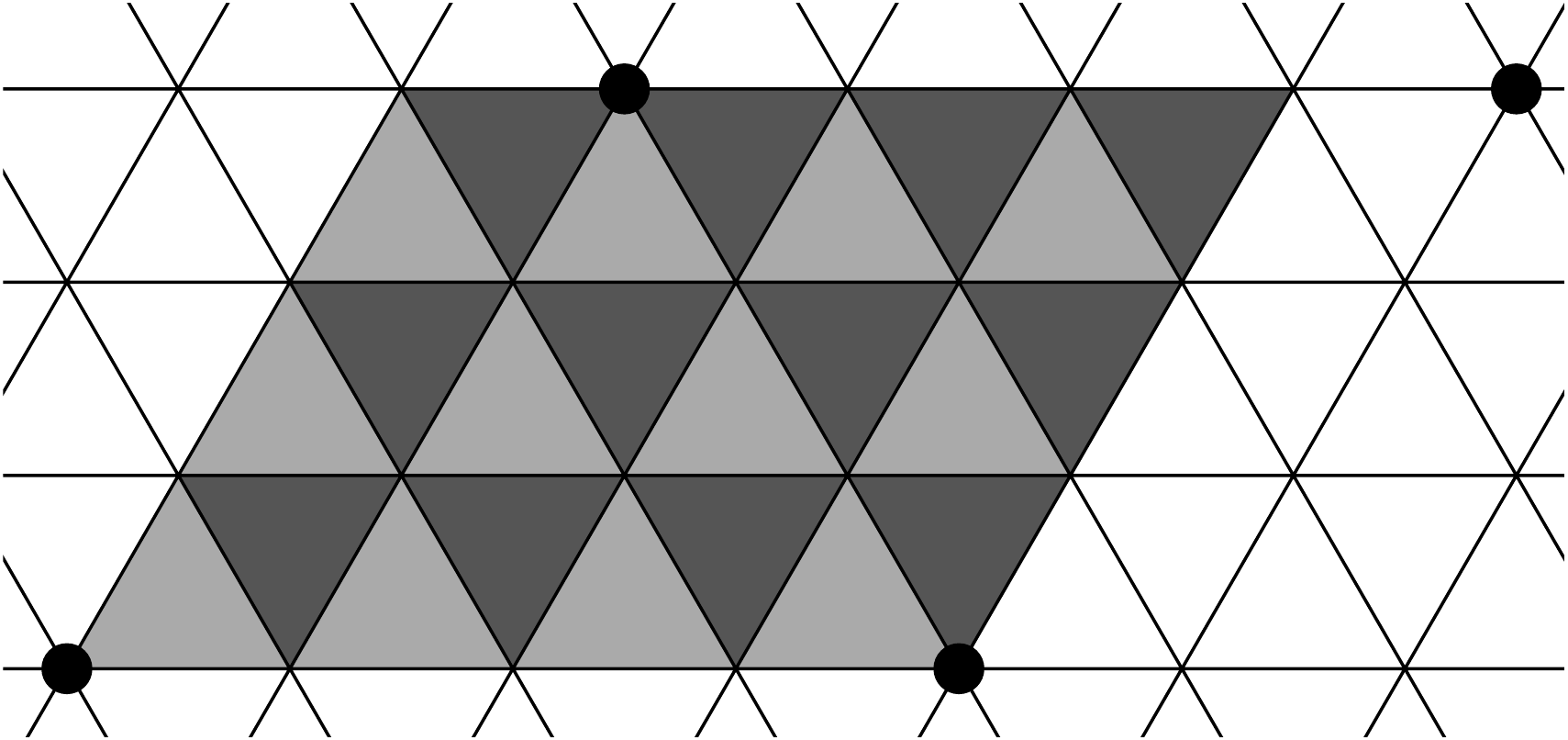}&\includegraphics[height=2.7cm]{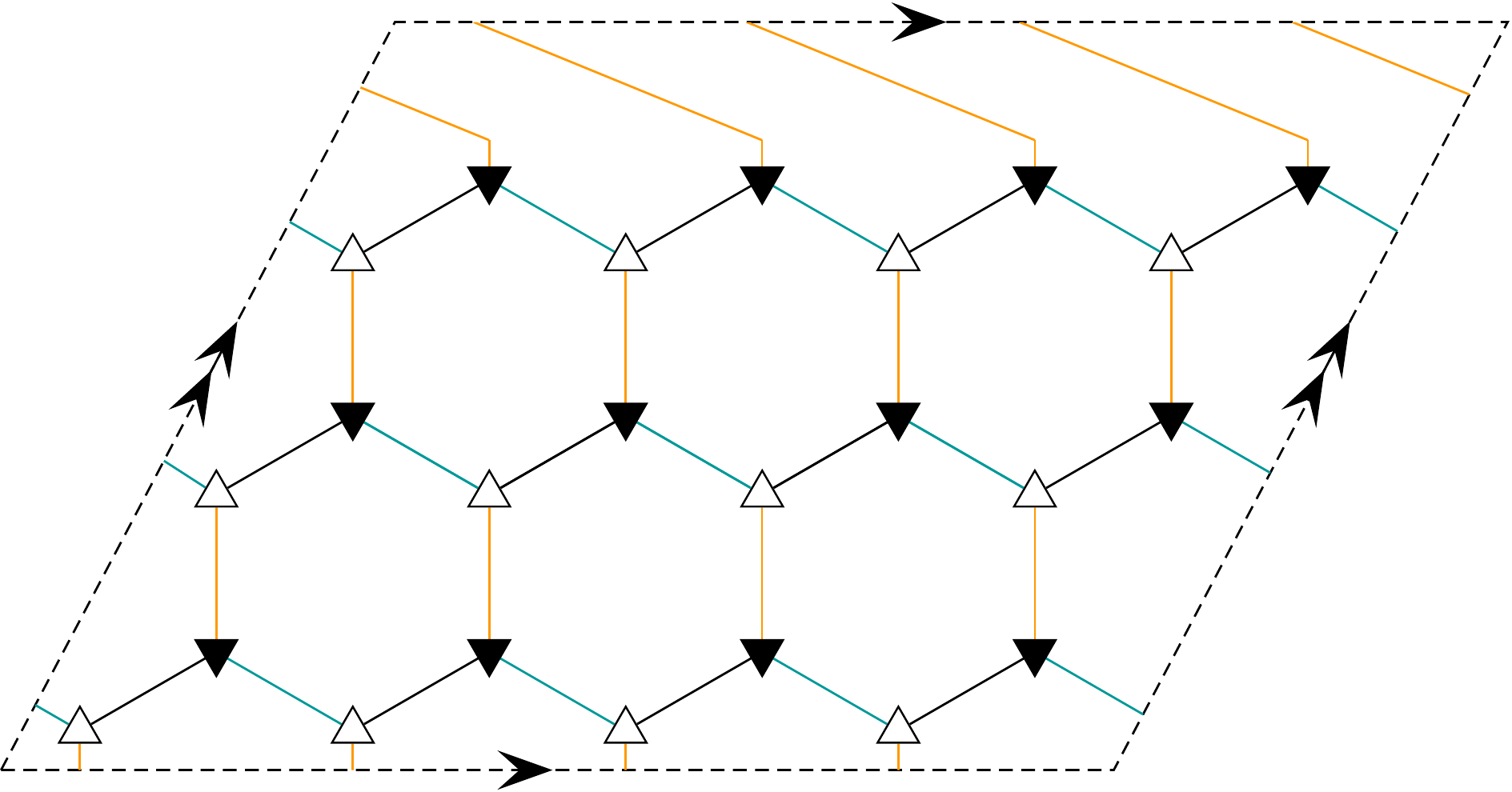}
\end{tabular}
\caption{The bipartite graph associated to a sublattice of \(\L_0\)}\label{eje_1}
\end{center}
\end{figure}

In this article we focus the problem of enumerating the
\(\L\)-periodic tilings, for a given full-rank sublattice of
\(\L_0\). Indeed, these tilings correspond to perfect matchings in a
certain ``honeycomb-like'' graph (see Figure~\ref{eje_1}), obtained as
follows. We consider two sets of representatives of \(\L_0/\L\) such
that the set of associated \up and \down triangles is connected. This
way, from \(\L\) we build (in principle not in a unique way) a graph
whose vertices are the triangles in a fundamental region:
\(V=\L_0/\L\times\L_0/\L\). An edge (weighted if desired by L,D,R)
joins a pair of triangles which can be merged in a lozenge.

The Hafnian of a symmetric square matrix of order \(2n\) is defined by:
\[\mathrm{Hf}(A)=\sum_{m\in U_n}\prod_{\{i,j\}\in m}A_{i,j},\]
where \(U_n=\left\{\{\{i_1,j_1\},\ldots,\{i_n,j_n\}\}\ |\ \cup_{1\leq
  l\leq n}\{i_l,j_l\}=\{1,\ldots,2n\}\right\}\) is the set of the
\((2n)!/(n!2^n)\) (unordered) matchings of the set
\(\{1,\ldots,2n\}\). The number of perfect matchings in an undirected
graph with an even number of vertices (the problem keeps no interest
if the number of vertices is odd) equals the Hafnian of its adjacency
matrix, but the computation of a Hafnian is unfortunately a
\#P-complete problem~\cite{Val_79}. The Pfaffian is a related, but
tractable, function defined over antisymmetric matrices:
\[\mathrm{Pf}(A')=\sum_{m\in U_n}(-1)^{\sigma(m)}\prod_{\genfrac{}{}{0pt}{}{\{i,j\}\in m}{i<j}}A'_{i,j},\]
where the {\it parity} \(\sigma(m)\) of a matching
\(m=\{\{i_1,j_1\},\ldots,\{i_n,j_n\}\}\), written in such a way that
\(i_l<j_l\), for \(1\leq l\leq n\), is that of the permutation
\[\left(\begin{array}{ccccccc}
1&2&3&4&\cdots&2n-1&2n\\ i_1&j_1&i_2&j_2&\cdots&i_n&jn\end{array}\right).\]
It can be shown that, for every antisymmetric matrix \(A'\) of even
order, we have \(\det(A')=\mathrm{Pf}(A')^2\). This result also holds
for matrices of odd order, defining their Pfaffian as zero.

Given a directed graph, the Pfaffian of its adjacency matrix \(A'\)
counts its perfect matchings, but some of them affected by a negative
sign, leading to an (in principle) meaningless sum. Kasteleyn proved
(see~\cite{Kas_63}) that every {\it planar} graph can be oriented in
such a way that \(\mathrm{Hf}(A)=\pm\mathrm{Pf}(A')\), counting
therefore the perfect matchings. This {\it Pfaffian orientation} can
be achieved requiring that in every face, the number of border edges
in clockwise direction is odd.

However, the graph of Figure~\ref{eje_1} is not planar. In general,
the graphs we are interested in can be embedded in a
torus. In~\cite{Kas_61}, the perfect matching enumeration problem is
solved for the rectangular lattice on a torus, by computing a linear
combination of four Pfaffians. In a more general fashion, it is stated
in \cite{Kas_67} that the problem can be solved with \(4^g\) Pfaffians
for any graph drawn in a surface of genus \(g\). This statement is
proved in \cite{Gal_Loe_99}; and independently in \cite{Tes_00}, which
contains a general method that uses \(2^{2-\chi}\) Pfaffians for a
graph embedded in a surface of Euler characteristic \(\chi\),
improving therefore the previous statement for non-orientable
surfaces. Tesler method~\cite{Tes_00} starts drawing the graph in a
surface represented by a polygon with pasted borders, distinguishing
between edges contained in the interior of the polygon and those
crossing its borders. If some of the interior edges form a cycle
enclosing the rest, a {\it crossing orientation} is given to the graph
according to the following rule (R4):
\begin{itemize}
\item The set of interior edges is given a Pfaffian orientation, as a
  planar graph.
\item Each of the remaining edges is oriented in such a way that any
  face it forms with the interior edges has an odd number of clockwise
  edges as well.
\end{itemize}

Before proceeding, let us remark that the graphs we consider are
bipartite (\up may only match \down, and conversely). For bipartite
graphs, the Permanent-Determinant variant, introduced
in~\cite{Per_69}, simplifies the Hafnian-Pfaffian method. Indeed, if
\(G\) and \(\vec G\) are, respectively, an undirected bipartite and a
directed bipartite graph, let \(A\) and \(A'\) denote their adjacency
matrices, and \(B\) and \(B'\) their bipartite adjacency matrices
(whose columns and rows represent black and white vertices,
respectively). Then, \(\mathrm{Hf}(A)=\mathrm{Per}(B)\) and
\(\mathrm{Pf}(A')=\det(B')\).

Let \(\L\) be a full-rank sublattice \(\L\) of \(\L_0\). Firstly, we
compute the (unique) matrix
\[B=\left[\begin{array}{cc}a&c\\0&b\end{array}\right]\in\Z^{2\times
    2}\] such that \(a>0\), \(0\leq c<b\), and \(\L\) is generated by
  the columns of \([\ve u|\ve v]B\) (see Equation~(\ref{Lambda})). We
  choose sets of \up and \down representatives formed in both cases by
  \(b\) rows of \(a\) triangles, in the way depicted in
  Figure~\ref{eje_1}. More explicitly:
\[\begin{array}{l@{\quad}l}
\bigtriangleup:&\{i\ve v+j\ve u\ |\ 0\leq i<b,\ 0\leq j<a\},\\
\bigtriangledown:&\{(i+1)\ve v+j\ve u\ |\ 0\leq i<b,\ 0\leq j<a\}.
\end{array}\]

\begin{figure}
\begin{center}
\begin{tabular}{cc}
\includegraphics[height=3cm]{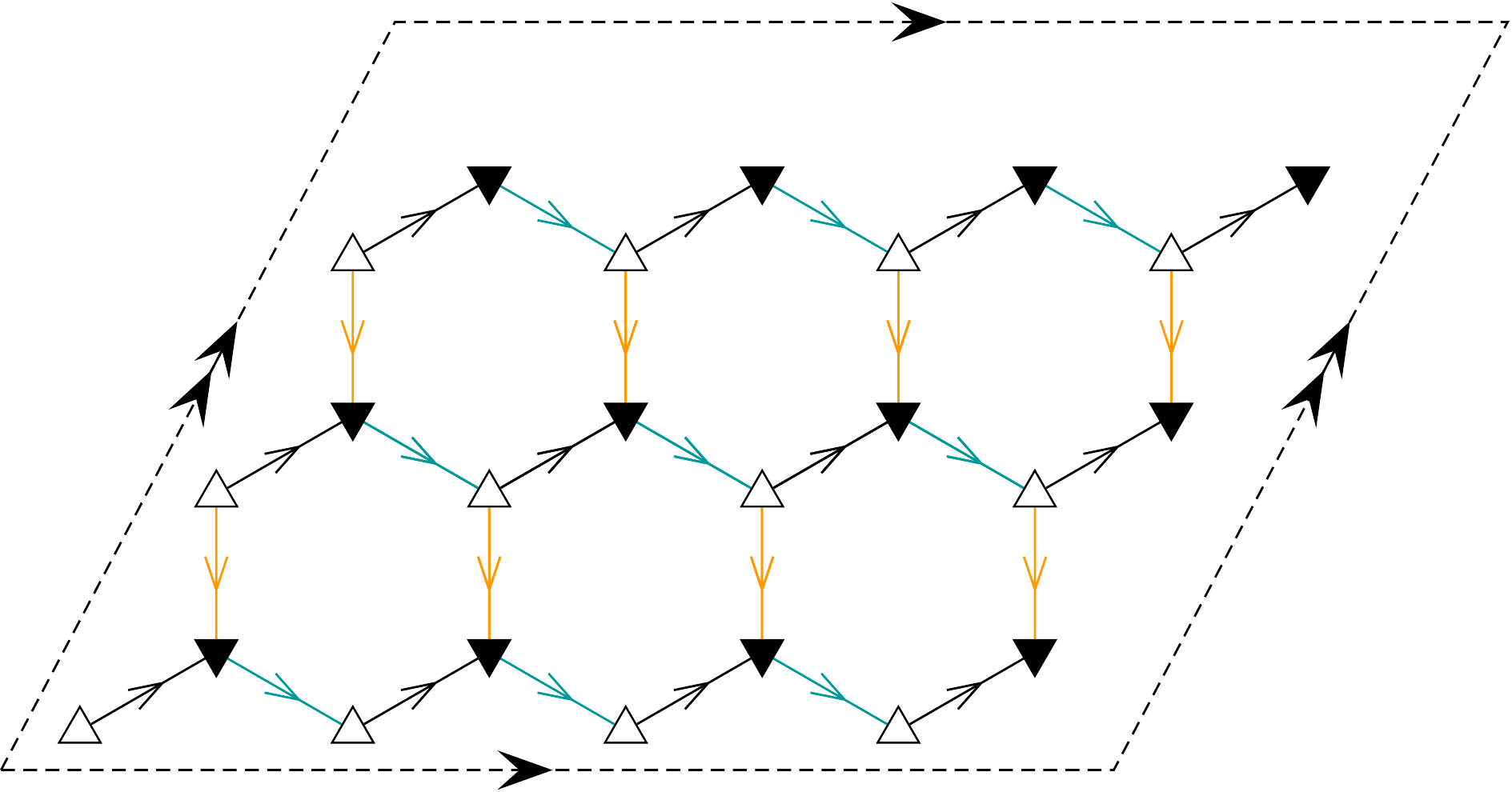}&\includegraphics[height=3cm]{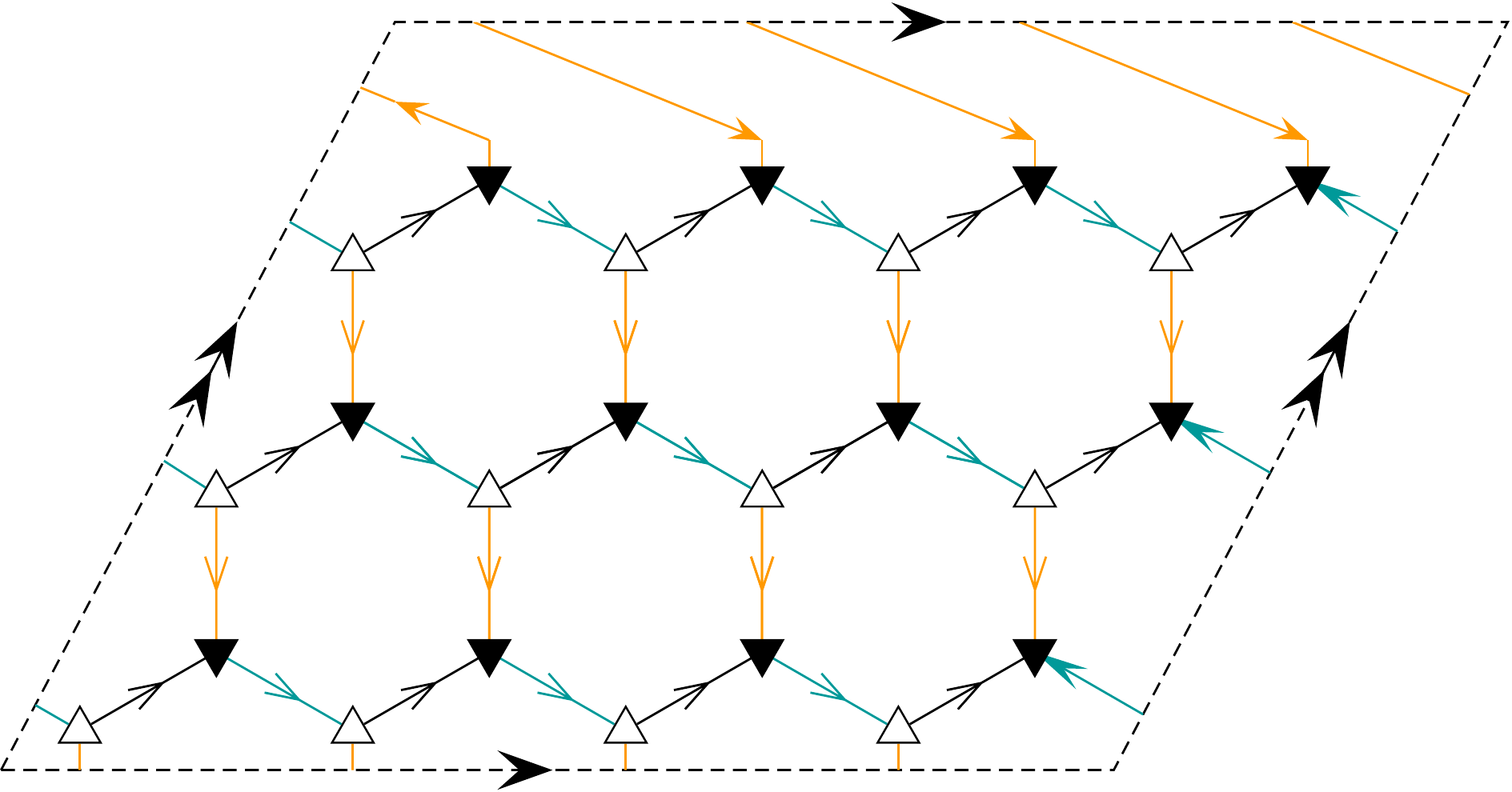}
\end{tabular}
\caption{Crossing orientation}\label{interiores}
\end{center}
\end{figure}

We label the edges of the graph with L,D, or R, depending of the
orientation of the corresponding lozenge. Then, the bipartite
adjacency matrix of the undirected graph, identifying rows with \up
and columns with \down, presents the following decomposition in
\(b\times b\) square blocks of order \(a\):
\[M=\left[\begin{array}{cccccc}
X&&&Z'\\
Z&X&&\\
&\ddots&\ddots&\\
&&Z&X\\
\end{array}\right],\]
where \(Z=\mathrm{D}\cdot\mathrm{Id}_a\),
\[X=\left[\begin{array}{cccc}
\mathrm{R}&&&\mathrm{L}\\
\mathrm{L}&\mathrm{R}&&\\
&\ddots&\ddots&\\
&&\mathrm{L}&\mathrm{R}\\
\end{array}\right],\quad\mathrm{and}\quad Z'=\left[\begin{array}{c|c}
&\mathrm{D}\cdot\mathrm{Id}_{a-c}\\ \hline
\mathrm{D}\cdot\mathrm{Id}_c&
\end{array}\right].\]
In the extreme case \(a=1\), block \(X\) equals \([R+L]\); and when
\(b=1\), the block decomposition of \(M\) collapses to
\([X+Z']\). Now, in order to define a suitable orientation, we start
by the edges which do not cross the border of the rectangle, orienting
by the rule \(\bigtriangleup\rightarrow\bigtriangledown\) the D and R
edges, and conversely the L edges, as is shown in the left side of
Figure~\ref{interiores}.

In order to apply Rule (R4) for orienting the rest of the edges, we
would need that those already oriented are enclosed in a cycle,
condition that is not fulfilled in our case. However, we can slightly
modify (R4) taking into account those D-edges which, joined to some
interior edges, form a cycle enclosing an odd number of vertices (no
R- or L-edges present this pathology under our construction). Then, we
get a crossing orientation if we orient the border crossing edges as
follows (see Figure~\ref{interiores}):
\begin{itemize}
\item L-edges: \(\bigtriangleup\rightarrow\bigtriangledown\).
\item D-edges:
\begin{itemize}
\item \(\mathrm{Id}_{a-c}\): \(\bigtriangleup\rightarrow\bigtriangledown\), if \(b\) is odd; and conversely otherwise.
\item \(\mathrm{Id}_{c\phantom{-a}}\): \(\bigtriangledown\rightarrow\bigtriangleup\), if \(b\) is odd; and conversely otherwise.
\end{itemize}
\end{itemize}
Marking edges crossing the diagonal border with \(\omega_1\) and those
crossing the horizontal one with \(\omega_2\), the block structure of
the bipartite adjacency matrix \(M'\) of this directed graph remains
as in \(M\), substituting the blocks by:
\(Z=\mathrm{D}\,\mathrm{Id}_a\);
\[X=\left[\begin{array}{cccc}
\mathrm{R}&&&\omega_1 \mathrm{L}\\
-\mathrm{L}&\mathrm{R}&&\\
&\ddots&\ddots&\\
&&-\mathrm{L}&\mathrm{R}\\
\end{array}\right],\ \mathrm{if}\ a>1,\ \mathrm{and}\ X=[R+\omega_1L]\ \mathrm{otherwise;}\]
\[Z'=\left[\begin{array}{c|c}
&(-1)^{b+1}\omega_2 \mathrm{D}\,\mathrm{Id}_{a-c}\\ \hline
(-1)^b\omega_1\omega_2 \mathrm{D}\,\mathrm{Id}_c&
\end{array}\right].\]
We obtain, as a corollary of~\cite[Theorem 5.2]{Tes_00}:
\begin{teo}
Let \(\L\) be a full-rank sublattice of \(\L_0\). If
\(g(\omega_1,\omega_2)=\det(M')\in\Z[L,D,R][\omega_1,\omega_2]\),
where \(M'\) is defined above, the generating function of the
\(\L\)-periodic tilings is:
\[Z(L,D,R)=\frac{1}{2}\left(g(1,1)+g(1,-1)+g(-1,1)-g(-1,-1)\right).\]
\end{teo}
It would be interesting to derive a ``closed formula'' for this
generating function, or at least for the number of \(\L\)-periodic
tilings \(Z(1,1,1)\).

\section{Different tiling types}\label{sec_tri}

The three constant mappings \(\t_L\), \(\t_D\), and \(\t_R\) are
\(\L\)-periodic tilings for any sublattice \(\L\subseteq\L_0\). The
identification with piles of cubes gives, for each of these tilings, a
plane (orthogonal to the \(Z\), \(X\), and \(Y\) axes respectively;
see Figure~\ref{edge} for the latter two).

\begin{figure}[h!]
\begin{center}
\includegraphics[width=6cm]{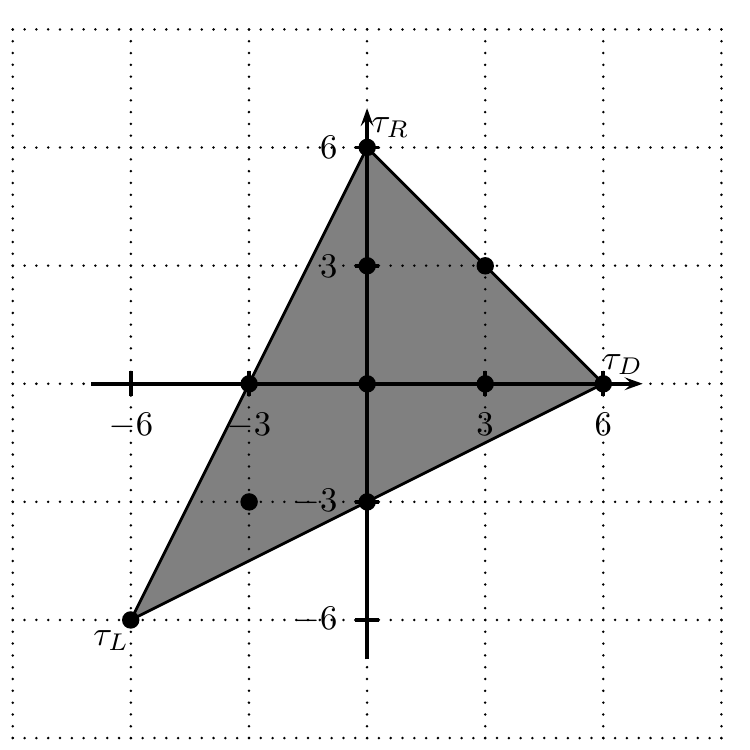}
\caption{Types of \(\L\)-periodic tilings.}\label{triangulo}
\end{center}
\end{figure}

Let \(B\in\Z^{2\times 2}\) be a 2-rank matrix. We call {\it
  fundamental triangle} to the triangle whose vertices are the
fingerprints in \(B\) of the constant tilings. With \(B=[\ve a|\ve
  b]\), we have:
\[\begin{array}{lcl}
\D(\t_L,B)&=&(-a_1-2a_2,-b_1-2b_2),\\
\D(\t_D,B)&=&(2a_1+a_2,2b_1+b_2),\\
\D(\t_R,B)&=&(-a_1+a_2,-b_1+b_2).
\end{array}\]
Therefore, the area of the fundamental triangle equals
\((9/2)\det\,B=(9/2)[\L_0:\L]\).  Note that the vertices of the
triangle are in the lattice \((3\Z)^2\), possibly shifted (see
Figure~\ref{triangulo}). Taking as example the lattice \(\L=[\ve u|\ve
  v]B\Z^2\), where
\[B=\left[\begin{array}{cc}2&-2\\ 2&4\end{array}\right],\]
the fundamental triangle is defined by the points:
\[\D(\t_L,B)=(-6,-6),\ \D(\t_D,B)=(6,0),\ \D(\t_R,B)=(0,6).\]

In general, for any two points \(\D_1,\D_2\in\Z^2\) such that
\(\D_1-\D_2\in(3\Z)^2\), there exists a third point
\(\D_3=-(\D_1+\D_2)\) defining with the former the fundamental
triangle associated to a basis.

Next result identifies the type of a tiling with a point of the
fundamental triangle. In the following, when speaking of points in the
fundamental triangle, we restrict ourselves to those in
\(\D(\t_D,B)+(3\Z)^2\).

\begin{teo}\label{teorema}Let \(\L\subseteq\L_0\) be a 2-rank lattice and \(\t\) a
  \(\L\)-periodic tiling. Writing \((L,D,R)=t(\t,\L)\) for the type
  of \(\t\), its fingerprint in a base \(B\) of \(\L\) has triangular
  coordinates proportional to \((L,D,R)\):
\[\D(\t,B)=\frac{1}{[\L_0:\L]}(L\D(\t_L,B)+D\D(\t_D,B)+R\D(\t_R,B)).\]
\end{teo}

\begin{demo}Let \(B=[\ve a|\ve b]\) be the considered basis. We denote by \(o_1\)
and \(o_2\) the orders in the group \(\L_0/\L\) of \(\ve u\) and
\(\ve v\), respectively, and consider the following two
representations of \(\L_0/\L\):
\[\begin{array}{|c|c|c|c|c|}\hline
0&\ve u&2\ve u&\cdots&(o_1-1)\ve u\\ \hline
\ve x_1&\ve x_1+\ve u&\ve x_1+2\ve u&\cdots&\ve x_1+(o_1-1)\ve u\\ \hline
\vdots&\vdots&\vdots&\ddots&\vdots\\ \hline
\ve x_{l_1-1}&\ve x_{l_1-1}+\ve u&\ve x_{l_1-1}+2\ve u&\cdots&\ve x_{l_1-1}+(o_1-1)\ve u\\ \hline
\end{array}\]
\[\begin{array}{|c|c|c|c|c|}\hline
0&\ve v&2\ve v&\cdots&(o_2-1)\ve v\\ \hline
\ve y_1&\ve y_1+\ve v&\ve y_1+2\ve v&\cdots&\ve y_1+(o_2-1)\ve v\\ \hline
\vdots&\vdots&\vdots&\ddots&\vdots\\ \hline
\ve y_{l_2-1}&\ve y_{l_2-1}+\ve v&\ve y_{l_2-1}+2\ve v&\cdots&\ve y_{l_2-1}+(o_2-1)\ve v\\ \hline
\end{array},\]
where \(l_i=[\L_0:\L]/o_i\). For \(i\in\{1,2\}, 0\leq j<l_i\) , let
\(L_j^i,D_j^i,\) and \(R_j^i\) be the number of elements in the
\(j\)th row of the \(i\)th table whose image by \(\tau\) is \(L,D,\) and \(R\), respectively. We have:
\[L_j^i+D_j^i+R_j^i=o_i,\forall i,\forall j.\]
There exist constants
\(\kappa_1=(\kappa_1^1,\kappa_1^2,\kappa_1^3)=H_\t(o_1\ve u),
\kappa_2=(\kappa_2^1,\kappa_2^2,\kappa_2^3)=H_\t(o_2\ve v)\) such
that
\[\kappa_1^1=-L_j^1-R_j^1,\ \kappa_1^2=\kappa_1^3=D_j^1,\]
\[\kappa_2^1=\kappa_2^2=-L_j^2,\ \kappa_2^3=D_j^2+R_j^2.\]
Therefore, the numbers \(D_j^1,L_j^2,L_j^1+R_j^1,\) and
\(D_j^2+R_j^2\) are independent of index \(j\). We consider the matrix
\[P=B^{-1}\left[\begin{array}{cc}o_1&\\ &o_2\end{array}\right]=
\left[\begin{array}{cc}b_2l_1&-b_1l_2\\-a_2l_1&a_1l_2\end{array}\right].\]
Writing \(k_i=\kappa_i^1+\kappa_i^2+\kappa_i^3\), we have
\((k_1,k_2)=(e_\t(\ve a),e_\t(\ve b))P\), and therefore,
\[\left\{\begin{array}{l}
e_\t(\ve a)=a_1k_1/o_1+a_2k_2/o_2\\
e_\t(\ve b)=b_1k_1/o_1+b_2k_2/o_2,
\end{array}
\right.\]
\[\begin{array}{ll}
{[\L_0:\L]}e_\t(\ve
a)=&\left(a_1l_1(2D_j^1-L_j^1-R_j^1)+a_2l_2(D_j^2-2L_j^2+R_j^2)\right)=\\
&a_1(2D-L-R)+a_2(D-2L+R),\\
{[\L_0:\L]}e_\t(\ve b)=&\left(b_1l_1(2D_j^1-L_j^1-R_j^1)+b_2l_2(D_j^2-2L_j^2+R_j^2)\right)=\\
&b_1(2D+L+R)+b_2(D-2L+R).\\
\end{array}\]
The result follows easily.
\end{demo}

For instance, the tiling shown in Figure~\ref{mixta} is
\(\L\)-periodic where \(\L\) is defined in the example from
Figure~\ref{triangulo}. Its type is (2,2,8); and its fingerprint in
basis \(B\), (0,3).

According to Figure~\ref{triangulo}, there are (at most) seven
nonconstant types of \(\L\)-periodic tilings. Four of them involve
lozenges in the three different orientations. Let us prove that every
point \(\D\) in the fundamental triangle represents at least one
tiling.

\begin{teo}\label{alpha}
Let \(B=[\ve a|\ve b]\) be a 2-rank integer matrix and \(\L=[\ve
u|\ve v]B\Z^2\). The set of types \(\{t(\t,\L)\ |\ \t\in T_\L\}\)
coincides with the intersection of the fundamental triangle with
\(\D(\t_D,B)+(3\Z)^2\).
\end{teo}
\begin{demo}
As a corollary of Theorem~\ref{teorema}, every type lies in that
intersection. For the converse inclusion, let \(\D\) be a point in the
fundamental triangle. According to Proposition~\ref{esqueleto}, the
``skeleton'' of any tiling with fingerprint \(\D\) is determined. We
claim that there is a tiling \(\t\) with the ``skeleton'' determined
by \(\D\). More formally, there is a \(\L\)-periodic tiling \(\t\)
such that the three-dimensional coordinates of the points \(a_1\ve
u+a_2\ve v\) and \(b_1\ve u+b_2\ve v\) are, respectively,
\begin{eqnarray*}\Xi(\t,B)=\left(\frac{1}{3}\left(\D_1-2a_1-a_2,\D_1+a_1-a_2,\D_1+a_1+2a_2\right),\right.\\
 \left.\frac{1}{3}\left(\D_2-2b_1-b_2,\D_2+b_1-b_2,\D_2+b_1+2b_2\right)\right).\end{eqnarray*}
A particular tiling can be formed just by placing a cube under every
point of that skeleton and filling the position whose coordinates are
not bigger component-wisely. For instance, applying this
process on the skeleton \(\langle(-2,0,2),(1,-1,3)\rangle\) we get the
tiling depicted in Figure~\ref{mixta}. We just need to show that for
every \(\ve x,\ve y\in\L\), \(H_\t(\ve x)-H_\t(\ve
y)\not\in(\N\backslash\{0\})^3\). In other words, that no point in the
skeleton is ``hidden''. It is sufficient to check that the normal to the
plane defined by the two components of \(\Xi(\tau,B)\) lies in the
cone \(\{\pm(x,y,z)\in\R^3\ |\ x\geq 0,y\geq 0,z\geq 0\}\). This
condition holds, because of Theorem~\ref{teorema} and the fact that
the cone is a convex set.
\end{demo}

The number of points in the border of the fundamental triangle are:
\[\begin{array}{cl}
\t_D,\t_L:&\gcd(a_1+a_2,b_1+b_2)+1\\
\t_L,\t_R:&\gcd(a_2,b_2)+1\\
\t_R,\t_D:&\gcd(a_1,b_1)+1\\
\end{array}\]

By Pick's Theorem, the number of interior points in the fundamental
triangle equals
\[\frac{1}{2}[\L_0:\L]-\frac{1}{2}(\gcd(a_1,b_1)+\gcd(a_2,b_2)+\gcd(a_1+a_2,b_1+b_2))+1\]
and the number of monomials in the generating function \(Z(L,D,R)\) is
\[\frac{1}{2}[\L_0:\L]+\frac{1}{2}(\gcd(a_1,b_1)+\gcd(a_2,b_2)+\gcd(a_1+a_2,b_1+b_2))+1.\]

\section{Grouping similar tilings}\label{sec_group}

The coefficients in the monomials corresponding to borders of the
triangle (i.e. types involving just one or two lozenge orientations)
are easily determined. Let us show the case of tiling types with no
``L'' lozenge. Put \(d=\gcd(a_1,b_1)\). There are \(d+1\) points in
the edge limited by \(\D(\t_D,B)\) and \(\D(\t_R,B)\). These are:
\[\D_i=\frac{i}{d}\D(\t_D,B)+\frac{d-i}{d}\D(\t_R,B),\
i=0,\ldots,d.\] According with Theorem~\ref{teorema}, a
\(\L\)-periodic tiling \(\t\) such that \(\D(\t,B)=\D_i\) has type
\((i[\L_0:\L]/d,0,(d-i)[\L_0:\L]/d)\). This kind of tilings are
constant on the lines directed by \(\ve v\) and can be enumerated as
follows:
\begin{itemize}
\item Compute a triangular form of \(B\):
\[\left[\begin{array}{cc}
d&0\\
\phantom{}*&*
\end{array}\right].\]

\item Select \(i\) elements \(\ve x\in\{0,\ve u,2\ve
  u,\ldots,(d-1)\ve u\}\), and define \(\t(\ve x)=D\) for them. For
  the rest, set \(\t(\ve x)=R\).
\end{itemize}
Therefore, the number of tilings of type
\((i[\L_0:\L]/d,0,(d-i)[\L_0:\L]/d)\) is \(\genfrac{(}{)}{0pt}{}d i\) and there
are \(2^d\) tilings whose type is in the considered edge.

\begin{figure}[h!]
\begin{center}
\begin{tabular}{cc}
\includegraphics[width=3cm]{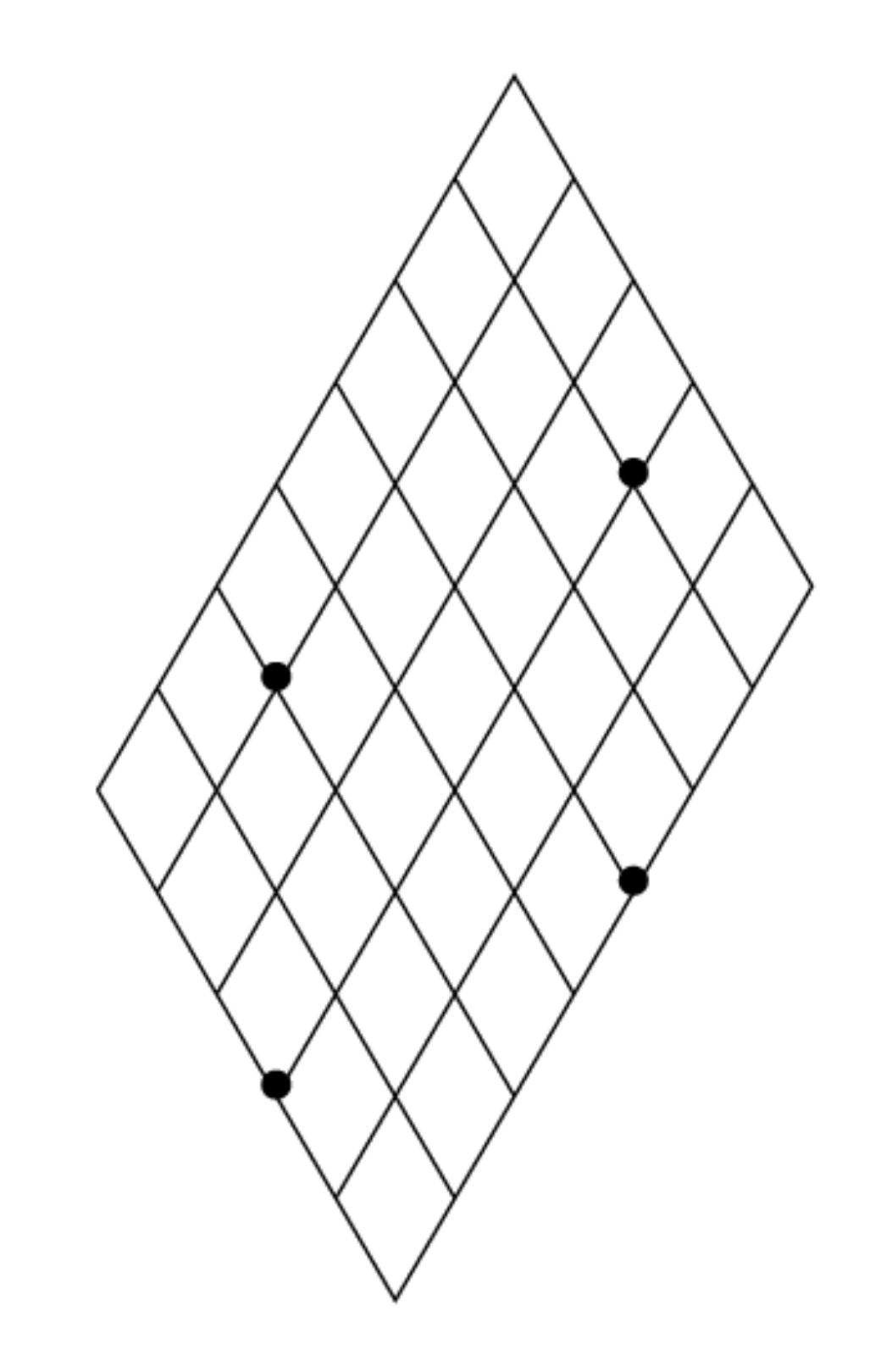}&\includegraphics[width=3cm]{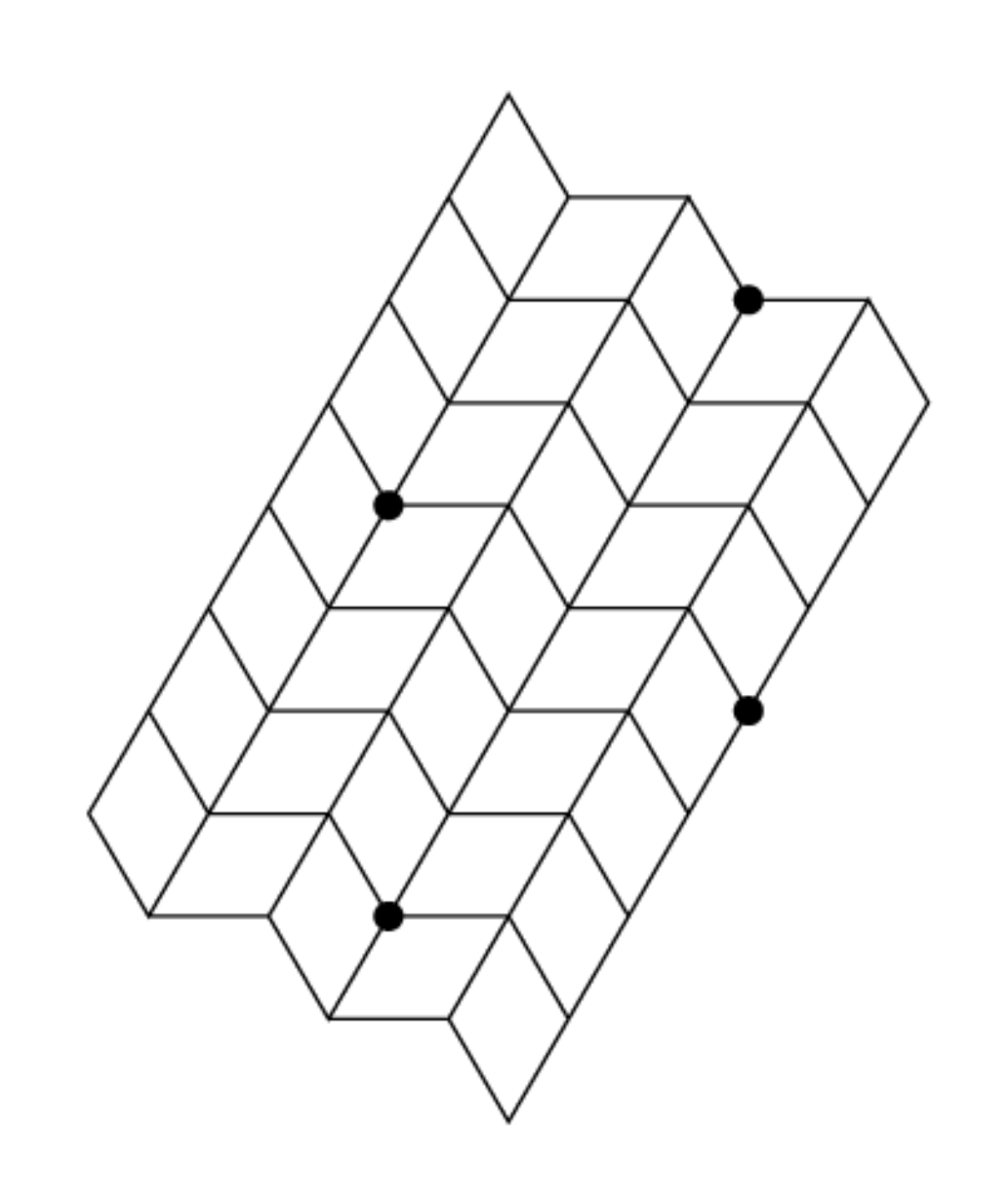}\\
\includegraphics[width=3cm]{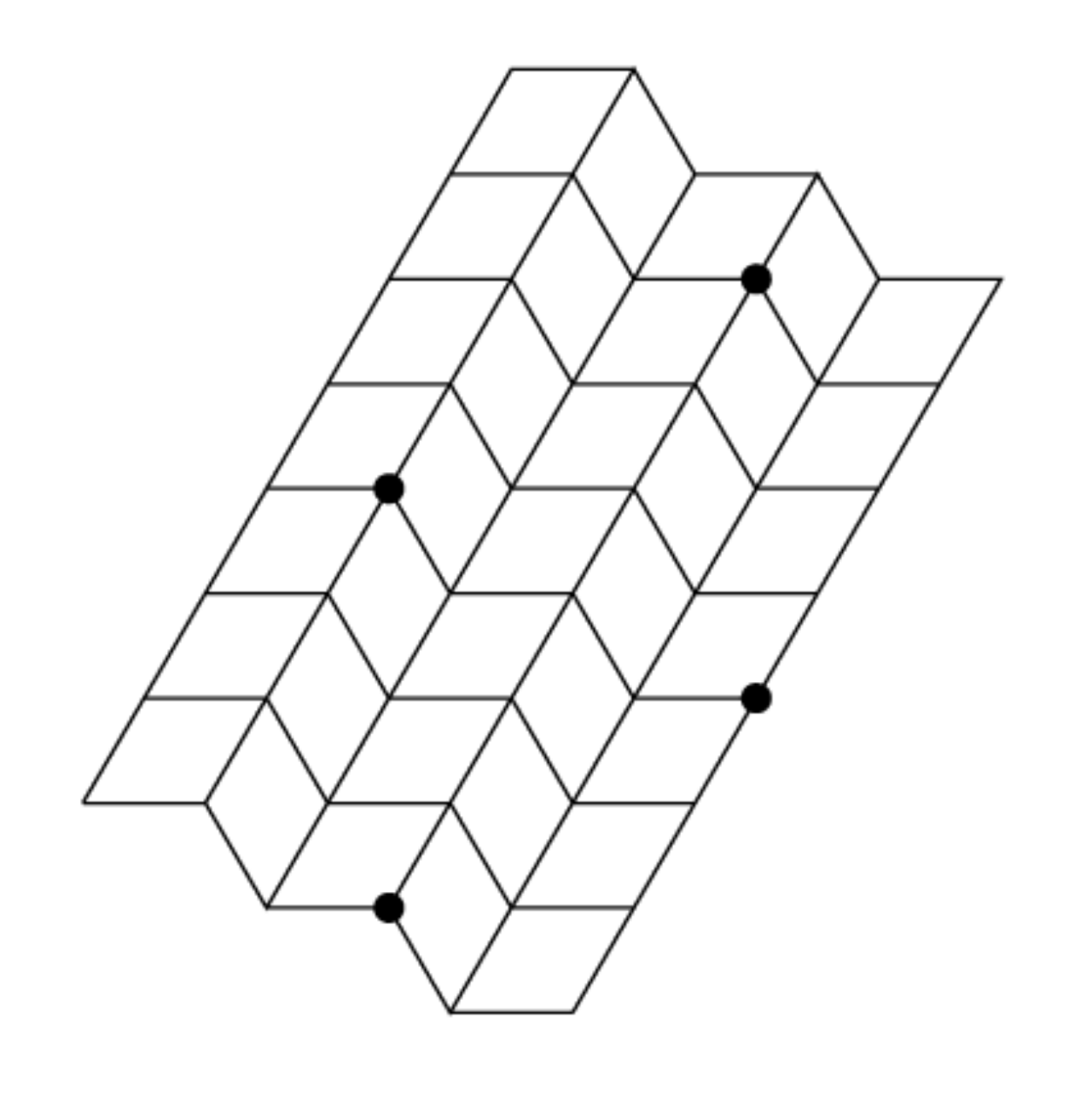}&\includegraphics[width=3cm]{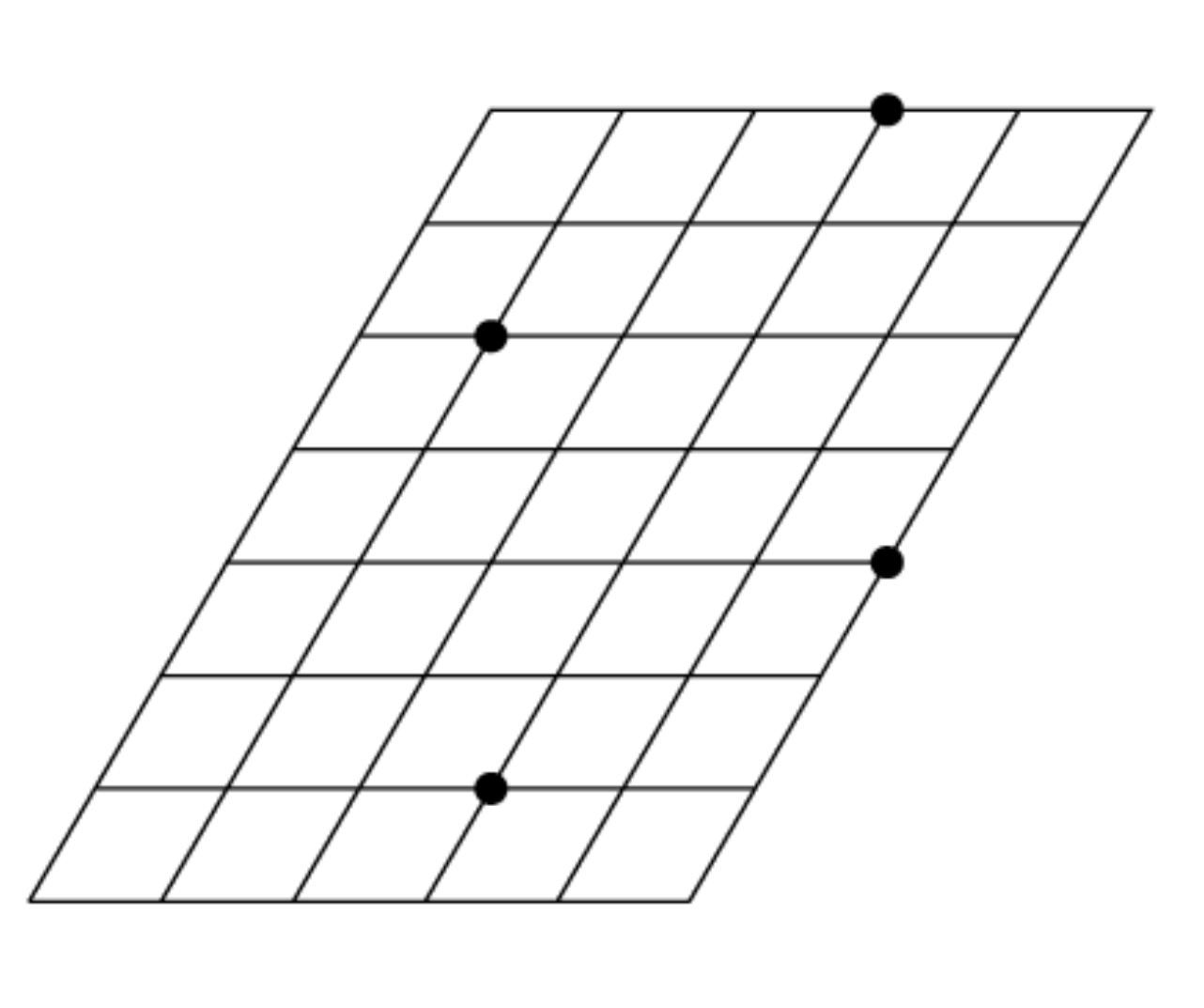}\\
\end{tabular}
\end{center}
\caption{Tilings with only ``D'' and ``R'' lozenges.}\label{edge}
\end{figure}

Continuing with the example from Figure~\ref{triangulo}, let us
enumerate the possible tilings in with no ``L'' lozenge. We need to
choose between \(D\) and \(R\) for \(\t(0)\) and \(\t(\ve u)\). We
get the four tilings depicted in Figure~\ref{edge}.

Indeed, the two tilings of type \(\D_1=(0,6,6)\) (with fingerprint
(3,3)) are the same modulo a shift. In general, we might also be
interested in enumerating classes of tilings, identifying those which
only differ on a shift. This is, defining the following equivalence
relation on the set of \(\L\)-periodic tilings \(T_\L\):
\[\tau S\rho\iff\exists\ve s\in\L_0\forall\ve x\in\L_0\ :\ \tau(\ve x)=\rho(\ve x+\ve
s),\] we need to count elements of \(T_\L/S\). In the degenerated case
of a triangle border, the number of classes associated to the point
\(\Delta_i\) is
\[\frac{1}{d}\sum_{k\mid(i,d)}\varphi(k)\genfrac{(}{)}{0pt}{}{d/k}{i/k},\]
the number of necklaces with \(d\) beans, \(i\) of them coloured. The
total number of classes of tilings modulo shifts in an edge is
\(\displaystyle\frac{1}{d}\sum_{k|d}\varphi(k)2^{d/k}\).

In this way, we can define a generating function
\(Z_1(L,D,R)\in\Z[L,D,R]\) of \(T_\L/S\), analogue to \(Z(L,D,R)\),
with the same set of monomials indeed.

Let us consider now an involution which associates pairs of tilings
with the same type. It is easily derived from
Definition~\ref{def_tiling} that a \(\L\)-periodic tiling can be
defined through a bijection \(\tilde\tau\) in \(\L_0\) such that
\(\tilde\tau(\ve x)-\ve x\in\{\ve 0,\ve v-\ve u,\ve v\}\) and
\(\tilde\tau(\ve x+\ve y)=\tilde\tau(\ve x)+\ve y\), for all
\(x\in\L_0\), \(y\in\L\). We define then the involution \(I\) as
follows:
\[(I\tilde\tau)(\ve x)=-\tilde\tau^{-1}(-\ve x).\]
This operation is compatible with the relation \(S\) defined above and
it is indeed more natural to consider \(I\) acting on \(T_\L/S\). In
Section~\ref{sec_def}, we set the convention that the left image in
Figure~\ref{esquinas} represents a solic cube. Considering the inverse
convention corresponds to looking at the pile of cubes ``from
behind''. This change of viewpoint is encoded by the involution
\(I\). Another interpretation arises from rotating \(180^\circ\) the
plane representation of the tiling.

We may also find redundant to compute as different tilings related by
this involution. This allows another simplification in the set
\(T_\L/S\), defining a coarser partition. In the degenerated case, the
cardinality of the new quotient is the number of reversible necklaces
with \(d\) beads, \(i\) of them coloured, but as before, the
computation of the generating function \(Z_2(L,D,R)\) seems a more
difficult problem.

As we have seen, for types lying in the border of the fundamental
triangle these functions correspond to well-known combinatorics
formulas. However, we are not able to efficiently compute them in the
more interesting case of points interior to the triangle.

Let \(\t\) be a \(\L\)-periodic tiling containing at least one lozenge
in each orientations (i.e., no component in its type is zero). This is
equivalent to the condition of having an inner corner (see
Figure~\ref{esquinas}) in the tiling. We call {\it flip} to a
transformation of one tiling into another by changing a inner corner
into a solic cube, or vice versa. It is clear that this operation
keeps the type of a tiling.

\begin{figure}[h!]
\begin{center}
\includegraphics[width=5cm]{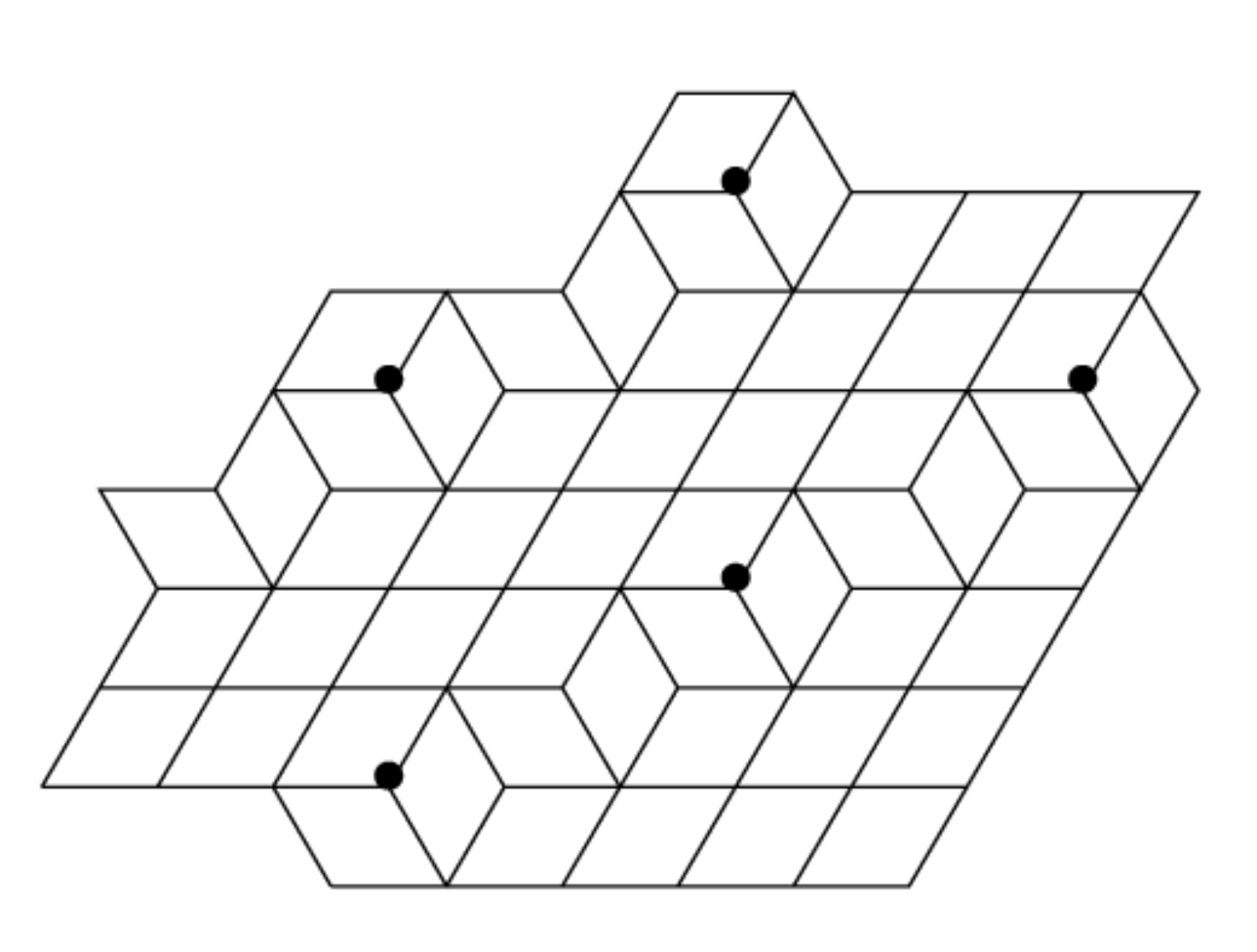}
\caption{\(\D(\tau,B)=(0,3)\)}\label{flipada}
\end{center}
\end{figure}

For example, a flip in the dotted quoin in the tiling from
Figure~\ref{mixta} (indeed, that tiling has only a quoin and an inner
corner) gives the tiling depicted in Figure~\ref{flipada}, which has
two quoins and two inner corners.

It follows from the proof of Theorem~\ref{alpha} that starting with a
\(\L\)-periodic tiling, one may obtain all the tilings with the same
type by means of flips. Perhaps this consideration is useful in the
task on determining \(Z_1,Z_2\). On the other hand, it is likely that
the Hafnian-Pfaffian method, which has proved useful in enumerating
symmetry classes of tilings (see~\cite{Kup_94}), facilitates the
evaluation of those functions.
\vspace{1cm}

\noindent {\bf Acknowledgements.} The author is partially supported by
the MTM2007-67088 research project of Spanish Ministerio de Ciencia e
Innovaci\'on. He thanks Professor Francisco Santos for his useful
comments on this article.

\end{document}